\documentclass[]{amsart}

\usepackage[initials]{amsrefs}
\usepackage{amsmath,amssymb}
\usepackage{amsthm}
\usepackage{combelow}
\usepackage{mathtools}
\usepackage{paralist}

\newtheorem{lem}{Lemma}
\newtheorem{prop}[lem]{Proposition}
\newtheorem{theo}[lem]{Theorem}
\newtheorem{cor}[lem]{Corollary}

\theoremstyle{definition}

\numberwithin{lem}{section}
\numberwithin{equation}{section}

\providecommand{\norm}[1]{\left\lVert#1\right\rVert}
\renewcommand{\Re}{\operatorname{Re}}

\begin{document}

\title{The almost Daugavet property and translation-invariant subspaces}
\author[S. L\"ucking]{Simon L\"ucking}
\address{Department of Mathematics\\ Freie Universit\"at Berlin\\ 
	 Arnimallee 6\\ 14195 Berlin\\ Germany}
\email{simon.luecking@fu-berlin.de}


\begin{abstract}
	Let $G$ be a metrizable, compact abelian group and let $\Lambda$ be a 
	subset of its dual group $\widehat G$. We show that $C_\Lambda(G)$ has the
	almost Daugavet property if and only if $\Lambda$ is an infinite set, and
	that $L^1_\Lambda(G)$ has the almost Daugavet property if and only if
	$\Lambda$ is not a $\Lambda(1)$ set.
\end{abstract}

\maketitle

\section{Introduction}

I. K. Daugavet \cite{Daugavet} proved in 1963 that all compact operators $T$
on $C[0,1]$ fulfill the norm identity
\begin{equation*}
	\norm{\mathrm{Id} + T} = 1 + \norm{T},
\end{equation*}
which has become known as the \emph{Daugavet equation}. C. Foia{\cb{s} and
I. Singer \cite{FoiasSingerPointsDiffusion} extended this result to all 
weakly compact operators on $C(K)$ where $K$ is a compact space without 
isolated points. Shortly afterwards G. Ya. Lozanovski{\u\i}
\cite{LozanovskiiAlmostIntegralOperators} showed that the Daugavet equation
holds for all compact operators on $L^1[0,1]$ and J. R. Holub
\cite{HolubDaugavetsEquationL1} extended this result to all weakly compact
operators on $L^1(\mu)$ where $\mu$ is a $\sigma$-finite non-atomic measure.
V. M. Kadets, R. V. Shvidkoy, G. G. Sirotkin, and D. Werner
\cite{KadetsShvidkoySirotkinWernerDaugavetProperty} proved that the validity
of the Daugavet equation for weakly compact operators already follows from the
corresponding statement for operators of rank 1. This result led to the
following definition: A Banach space $X$ is said to have the \emph{Daugavet
property}, if every operator $T: X \rightarrow X$ of rank 1 satisfies the
Daugavet equation. During the studies of ultraproducts
\cite{BilikKadetsShvidkoyWernerNarrowOperatorsDaugavetPrUltraproducts}
and quotients \cite{KadetsShepelskaWernerQuotientsDaugavetProperty} of Banach
spaces with the Daugavet property a weaker notion was introduced. Let $X$ be a 
Banach space and let $Y$ be a subspace of $X^*$. We say that $X$ has the
\emph{Daugavet property with respect to $Y$}, if the Daugavet equation holds 
true for every rank-one operator $T: X \rightarrow X$ of the form 
$T=y^* \otimes x$ where $x \in X$ and $y^* \in Y$. A Banach space $X$ is 
called an \emph{almost Daugavet space} or a space with the \emph{almost
Daugavet property}, if it has the Daugavet property with respect to some 
norming subspace $Y \subset X^*$. Recall that a subspace $Y \subset X^*$ is
said to be norming, if for every $x \in X$
\begin{equation*}
	\sup_{y^* \in S_Y} \left| y^*(x) \right| = \norm{x}.
\end{equation*}
The space $\ell^1$ is an example of an almost Daugavet space that does not 
have the Daugavet property.

Separable almost Daugavet spaces can be characterized using a kind of inner
measure of non-compactness of the unit sphere. We call a set $F$ an inner
$\varepsilon$-net for $S_X$, if $F \subset S_X$ and for every $x \in S_X$ there
is a $y \in F$ with $\norm{x - y} \leq \varepsilon$. Then the \emph{thickness}
$T(X)$ of a Banach space $X$ is defined by
\begin{equation*}
	T(X) = \inf \left\{ \varepsilon>0 : \text{there exists a finite
	inner $\varepsilon$-net for $S_X$} \right\}.
\end{equation*}
R. Whitley \cite{WhitleySizeUnitSphere} introduced this parameter and showed
that $1 \leq T(X) \leq 2$ if $X$ is infinite-dimensional, in particular that
$T(\ell^p) = 2^{1/p}$ for $1 \leq p < \infty$ and $T(C(K)) = 2$, if $K$ has no
isolated points. It was shown by V. Kadets, V. Shepelska, and D. Werner that
a separable Banach space $X$ is an almost Daugavet space if and only if 
$T(X)=2$ \cite{KadetsShepelskaWernerThicknessAlmostDaugavet}*{Theorem 1.1}.
 
Almost Daugavet spaces contain $\ell^1$
\cite{KadetsShepelskaWernerThicknessAlmostDaugavet}*{Corollary 3.3} and are
considered ``big''. It is therefore an interesting question which subspaces of
an almost Daugavet space inherit the almost Daugavet property. The most general
result in this direction is that a closed subspace $Z$ of a separable almost
Daugavet space $X$ has the almost Daugavet property as well, if the quotient
space $X/Z$ contains no copy of $\ell^1$ 
\cite{LueckingSubspacesAlmostDaugavet}*{Theorem 2.5}.

Let us consider an infinite, compact abelian group $G$ with its Haar measure
$m$. Since $G$ has no isolated points and since $m$ has no atoms, the spaces
$C(G)$ and $L^1(G)$ have the Daugavet property. Using the group structure of
$G$, we can translate functions that are defined on $G$ and look at closed,
translation-invariant subspaces of $C(G)$ or $L^1(G)$. These subspaces can be
described via subsets $\Lambda$ of the dual group $\widehat G$ and are of the
form
\begin{align*}
	C_\Lambda(G) &= \left\{ f \in C(G): \text{spec} f \subset \Lambda
		\right\}\\
	L^1_\Lambda(G) &= \left\{ f \in L^1(G): \text{spec} f \subset \Lambda
		\right\},
	\shortintertext{where}
	\text{spec}f &= \left\{ \gamma \in \widehat G: \hat f(\gamma) \neq
		0 \right\}.
\end{align*}
We are going to characterize the sets $\Lambda \subset \widehat G$ such that 
$C_\Lambda(G)$ and $L^1_\Lambda(G)$ are of thickness 2. If $G$ is metrizable,
this leads to a characterization of the translation-invariant subspaces of 
$C(G)$ and $L^1(G)$ which have the almost Daugavet property.

\section{Translation-invariant subspaces of $C(G)$}

Let us start with translation-invariant subspaces of $C(G)$. We will show that
$T(C_\Lambda(G)) = 2$ if and only if $\Lambda$ is an infinite subset of
$\widehat G$. We will split the proof into various cases that depend on the
structure of $G$. For this reason we recall some definitions and results
concerning abelian groups and compact abelian groups.

Let $G$ be an abelian group with identity element $e_G$. A subset $E$ of $G$ is
said to be \emph{independent}, if $x_1^{k_1}  \cdots  x_n^{k_n} =e_G$ implies
$x_1^{k_1} = \cdots = x_n^{k_n} = e_G$ for every choice of distinct points
$x_1, \dotsc, x_n \in E$ and integers $k_1, \dotsc, k_n$. The \emph{order}
$o(x)$ of an element $x \in G$ is the smallest positive integer $m$ such that
$x^m = e_G$. If no such $m$ exists, $x$ is said to have infinite order.

Let $\mathbb T$ be the \emph{circle group}, i.e., the multiplicative group of 
all complex numbers with absolute value one. If $G$ is a compact abelian group,
we denote the identity element of $\widehat G$, which coincides with the
function identically equal to one, by $\mathbf{1}_{G}$. Linear combinations of
elements of $\widehat G$ are called \emph{trigonometric polynomials} and for
every $\Lambda \subset \widehat G$ the space
$T_\Lambda(G) = \operatorname{lin} \Lambda$ of trigonometric polynomials with
spectrum contained in $\Lambda$ is dense in $C_\Lambda(G)$.

Let $H$ be a closed subgroup of $G$. The \emph{annihilator} of $H$ is defined
by
\begin{equation*}
	H^\perp = \left\{ \gamma \in \widehat{G}: \gamma(x) = 1 \text{ for all $x
	\in H$}	\right\}
\end{equation*}
and is therefore a closed subgroup of $\widehat G$. We have that
$\widehat H = \widehat G/H^\perp$ and that $\widehat{G/H} = H^\perp$
\cite{RudinFourierAnalysis}*{Theorem 2.1.2}.

If $(G_j)_{j \in J}$ is a family of abelian groups, we define their
\emph{direct product} (or their \emph{complete direct sum}) by
\begin{align*}
	\prod_{j \in J} G_j &= \left\{ f : J \rightarrow \bigcup_{j \in J} G_j:
		f(j) \in G_j  \right\}
	\shortintertext{and define the group operation coordinatewise. Their
		\emph{direct sum} is the subgroup}
	\bigoplus_{j \in J} G_j &= \left\{ f \in \prod_{j \in J} G_j :
		f(j) = e_{G_j} \text{ for all but finitely many $j \in J$} \right\}.
\end{align*}
If all $G_j$ coincide with $G$, we write $G^J$ or $G^{(J)}$ for the direct
product or the direct sum. We denote by $p_{G_j}$ the projection from
$\prod_{j \in J} G_j$ onto $G_j$. If we consider products of the form $\mathbb
Z^\mathbb N$ or $\mathbb Z^n$, we denote by $p_1, p_2, \dotsc$ the
corresponding projections onto $\mathbb Z$. If all $G_j$ are compact, then
$\prod_{j \in J} G_j$ is a compact abelian group as well and its dual group is
given by $\bigoplus_{j \in J} \widehat{G_j}$
\cite{RudinFourierAnalysis}*{Theorem 2.2.3}.
\begin{prop}
	\label{finiteTorus}
	Let $A$ be a compact abelian group, set $G = \mathbb T \oplus A$, and let 
	$\Lambda$ be a subset of $\widehat G = \mathbb Z \oplus \widehat A$. If
	$p_\mathbb Z[\Lambda]$ is infinite, then $T(C_\Lambda(G)) = 2$.
	\begin{proof}
		Fix $f_1, \dotsc, f_n \in S_{C_\Lambda(G)}$ and $\varepsilon > 0$. 
		We have to
		find $g \in S_{C_\Lambda(G)}$ with $\norm{f_k + g}_\infty \geq 2 -
		\varepsilon$ for $k = 1, \dotsc, n$.
		
		Every $f_k$ is uniformly continuous and therefore 
		there exists $\delta >0$ such that for $k = 1, \dotsc, n$ and 
		all $a \in A$
		\begin{equation*}
			|\varphi - \vartheta| \leq \delta \Longrightarrow
			\left| f_k(\mathrm e^{\mathrm i \varphi}, a) -
			f_k(\mathrm e^{\mathrm i \vartheta}, a) \right| \leq \varepsilon.
		\end{equation*}
		Since $p_\mathbb Z[\Lambda]$ contains infinitely many elements, we can
		pick $s \in	p_\mathbb Z[\Lambda]$ with $|s|2 \delta \geq 2\pi$. By our
		choice of $s$, we get for all $\vartheta \in [0, 2\pi]$
		\begin{equation}
			\label{finiteTorusEquation1}
			\left\{ \mathrm e^{\mathrm i s \varphi}:
			|\varphi - \vartheta| \leq \delta \right\}
			= \left\{ \mathrm e^{\mathrm i \varphi}: 
			|\varphi - \vartheta| \leq |s| \delta \right\} 
			= \mathbb T.
		\end{equation}
		Choose $g \in \Lambda$ with $p_\mathbb Z(g) = s$ and fix $k \in
		\{1,\dotsc, n\}$. There exists
		$(\mathrm e^{\mathrm i \vartheta^{(k)}}, a^{(k)}) \in G$ with
		\begin{equation*}
			\left| f_k(\mathrm e^{\mathrm i \vartheta^{(k)}}, a^{(k)}) \right|
			= 1,
		\end{equation*}
		since $f_k \in S_{C_\Lambda(G)}$. By 
		(\ref{finiteTorusEquation1}), we can pick $\varphi^{(k)} \in 
		\mathbb R$ with
		\begin{gather*}
		\left| \varphi^{(k)} - \vartheta^{(k)} \right| \leq \delta\\
		\shortintertext{and}
		\mathrm e^{\mathrm i s \varphi^{(k)}} =
		\frac{f_k(\mathrm e^{\mathrm i \vartheta^{(k)}}, a^{(k)})}
		{g(1,a^{(k)})}.
		\end{gather*}
		Note that the right-hand side of the last equation has absolute
		value 1 because $g$ is a character of $G$. Consequently,
		\begin{equation*}
			g(\mathrm e^{\mathrm i \varphi^{(k)}}, a^{(k)}) =
			g(\mathrm e^{\mathrm i \varphi^{(k)}}, e_A) g(1, a^{(k)}) =
			e^{\mathrm i s \varphi^{(k)}} g(1, a^{(k)}) =
			f_k(\mathrm e^{\mathrm i \vartheta^{(k)}}, a^{(k)}).
		\end{equation*}
		Finally,
		\begin{align*}
			\norm{f_k + g}_\infty &\geq
				\left| f_k(\mathrm e^{\mathrm i \varphi^{(k)}}, a^{(k)}) +
				g(\mathrm e^{\mathrm i \varphi^{(k)}}, a^{(k)}) \right|\\
			&\geq \left| g(\mathrm e^{\mathrm i \varphi^{(k)}}, a^{(k)}) +
				f_k(\mathrm e^{\mathrm i \vartheta^{(k)}}, a^{(k)}) \right|
				- \left| f_k(\mathrm e^{\mathrm i \varphi^{(k)}}, a^{(k)}) -
				f_k(\mathrm e^{\mathrm i \vartheta^{(k)}}, a^{(k)}) \right|\\
			&\geq 2 - \varepsilon. \qedhere
		\end{align*}
	\end{proof}	
\end{prop}
\begin{prop}
	\label{infiniteTorus}
	Let $A$ be a compact abelian group, set $G = \mathbb T^\mathbb N \oplus
	A$, and let $\Lambda$ be a subset of $\widehat G = \mathbb Z^{(\mathbb N)}
	\oplus \widehat A$. If we find arbitrarily large $l \in \mathbb N$ with
	$p_l[\Lambda] \neq \{0\}$, then $T(C_\Lambda(G)) = 2$.
	\begin{proof}
		Fix $f_1, \dotsc, f_n \in S_{C_\Lambda(G)}$. Since $T_\Lambda(G)$ is
		dense in $C_\Lambda(G)$, we may assume without loss of 
		generality that $f_1, \dotsc, f_n$ are trigonometric polynomials. 
		We are going to find $g \in S_{C_\Lambda(G)}$ with $\norm{f_k +
		g}_\infty = 2$ for $k = 1,\dotsc, n$.
		
		Setting $\Delta = \bigcup_{k=1}^n \text{spec}f_k$, we get a finite
		subset of $\Lambda$ because every $f_k$ is a trigonometric
		polynomial and therefore has a finite spectrum. Consequently, there
		exists $l_0 \in \mathbb N$ with $p_l[\Delta] = \{0\}$ for all
		$l > l_0$ and the evaluation of $f_1, \dotsc, f_n$ at a point
		$(t_1,t_2, \dotsc, a) \in G$ just depends on the coordinates
		$t_1, \dotsc, t_{l_0}$ and $a$.
		
		By assumption, we can find $l > l_0$ and $g \in \Lambda$ with
		$s = p_l(g) \neq 0$. Fix $k \in \{1,\dotsc, n\}$. There exists
		$x^{(k)} = (t^{(k)}_1,t^{(k)}_2, \dotsc, a^{(k)}) \in G$ with
		$|f_k(x^{(k)})|=1$ since $f_k \in S_{C(G)}$. Pick $u^{(k)} \in 
		\mathbb T$ 	with
		\begin{equation*}
			(u^{(k)})^s = \frac{f_k(x^{(k)})}
			{g(t^{(k)}_1, \dotsc, t^{(k)}_{l-1}, 1, t^{(k)}_{l+1},
			t^{(k)}_{l+2}, \dotsc, a^{(k)})}.
		\end{equation*}
		Note that the right-hand side of the last equation has absolute
		value 1 because $g$ is a character of $G$. With the same reasoning
		as at the end of the proof of Proposition \ref{finiteTorus} we
		get that
		\begin{equation*}
			g(t^{(k)}_1, \dotsc, t^{(k)}_{l-1}, u^{(k)}, t^{(k)}_{l+1},
			t^{(k)}_{l+2}, \dotsc, a^{(k)})=f_k(x^{(k)}).
		\end{equation*}
		Finally,
		\begin{align*}
			\norm{f_k + g}_\infty &\geq
				\left|(f_k + g)(t^{(k)}_1, \dotsc, t^{(k)}_{l-1}, u^{(k)},
				t^{(k)}_{l+1}, t^{(k)}_{l+2}, \dotsc, a^{(k)}) \right|\\
			&= 2 \left| f_k(x^{(k)}) \right| = 2. \qedhere
		\end{align*}
	\end{proof}
\end{prop}
\begin{lem}
	\label{Lemma1}
	Let $\varepsilon > 0$ and $z_1, \dotsc, z_n \in \{ z \in \mathbb C :
	|z| \leq 1 \}$ with
	\begin{equation*}
		\left| \sum_{k = 1}^n z_k \right| \geq n(1 - \varepsilon).
	\end{equation*}
	Then
	\begin{equation*}
		|z_k| \geq 1 - n\varepsilon	\quad \text{and} \quad
		|z_k - z_l| \leq 2n \sqrt \varepsilon	
	\end{equation*}
	for $k, l = 1, \dotsc, n$.
	\begin{proof}
		The first assertion is an easy consequence of the triangle
		inequality.
		
		For fixed $k, l \in \{1, \dotsc, n\}$ we have
		\begin{align*}
			\Re z_k \overline{z_l} &= \Re \sum_{\mathclap{s,t = 1}}^n 
				z_s \overline{z_t} - \Re \sum_{\mathclap{\substack{s,t=1\\
				(s,t) \neq (k,l)}}}^n z_s \overline{z_t} =
				\left| \sum_{k = 1}^n z_k \right|^{2} -
				\Re \sum_{\mathclap{\substack{s,t = 1\\(s,t) \neq (k,l)}}}^n
				z_s	\overline{z_t}\\
			&\geq n^2(1 - \varepsilon)^2 - (n^2 - 1) =
				1 - 2n^2\varepsilon + n^2\varepsilon^2\\
			&\geq 1 - 2n^2\varepsilon.
		\end{align*}
		Using this inequality, we get
		\begin{align*}
			|z_k - z_l|^2 &= |z_k|^2 + |z_l|^2 -2\Re z_k \overline{z_l}\\
			&\leq 2 - 2(1 - 2n^2\varepsilon) = 4n^2\varepsilon. \qedhere
		\end{align*}
	\end{proof}	
\end{lem}
The following lemma shows that if we are given $n$ subsets of the unit circle that do not
meet a circular sector with central angle bigger than $\frac{2\pi}{n}$, 
then we can rotate these $n$ subsets such that their intersection becomes empty.
\begin{lem}
	\label{Lemma2}
	Let $W_1, \dotsc, W_n \subset \{ z \in \mathbb C : |z| \leq 1 \}$. 
	Suppose that for every $k \in \{1, \dotsc, n\}$ there
	exist $\varphi_k \in [0, 2\pi]$ and $\vartheta_k \in [\frac{2\pi}{n},
	2\pi]$ with
	\begin{equation*}
		W_k \cap \{ r \mathrm e^{\mathrm i \alpha}: r \in [0,1],
		\alpha \in [\varphi_k, \varphi_k + \vartheta_k] \} = \emptyset.
	\end{equation*}
	Then there exist $t_1, \dotsc, t_n \in \mathbb T$ with
	\begin{equation*}
		\bigcap_{k = 1}^n t_k W_k = \emptyset.
	\end{equation*}
	\begin{proof}
		Setting for $k = 1, \dotsc, n$ (with $\vartheta_0 = 0$)
		\begin{gather*}
			t_k = \mathrm e^{\mathrm i \sum_{l = 0}^{k - 1}\vartheta_l}
			\mathrm e^{-\mathrm i \varphi_k},
			\shortintertext{we get}
			t_k W_k \cap
			\left\{ r \mathrm e^{\mathrm i \alpha}:
			r \in [0,1], \alpha \in
			\left[ \sum_{l = 0}^{k - 1}\vartheta_l,
			\sum_{l = 0}^{k}\vartheta_l \right] \right\} = \emptyset.
		\end{gather*}
		Fix $\alpha \in [0,2\pi]$ and $r \in [0,1]$. 
		Since $\sum_{k = 1}^n \vartheta_k \geq 2\pi$,
		there exists $k \in \{1, \dotsc, n\}$ with
		\begin{equation*}
			\alpha \in \left[ \sum_{l = 0}^{k - 1}\vartheta_l,
			\sum_{l = 0}^{k}\vartheta_l \right].
		\end{equation*}
		Consequently, $r \mathrm e^{\mathrm i \alpha}$ does not belong to
		$t_k W_k$ and $\bigcap_{k = 1}^n t_k W_k = \emptyset$.
	\end{proof}
\end{lem}
\begin{lem}
	\label{Lemma3}
	Let $\varepsilon, \delta > 0$, let $W \subset \{ z \in \mathbb C:
	1 - \delta \leq |z| \leq 1 \}$, and set $W_\varepsilon:
	\{ z \in \mathbb C: \text{there exists $w \in W$ with
	$|w-z| \leq \varepsilon$} \}$. Suppose that there exists $\vartheta
	\in [0,2\pi]$ such that for every $\varphi \in [0,2\pi]$
	\begin{equation*}
		W_\varepsilon \cap
		\{ r \mathrm e^{\mathrm i \alpha}: r \in [0,1], \alpha \in
		[\varphi, \varphi + \vartheta] \} \neq \emptyset.
	\end{equation*}
	Then $W$ is a $(2\varepsilon + \delta + \vartheta)$-net for $\mathbb T$.
	\begin{proof}
		Fix $\mathrm e^{\mathrm i \varphi} \in \mathbb T$. We have to find
		$w \in W$ with $|w- e^{\mathrm i \varphi}| \leq 2\varepsilon + \delta
		+ \vartheta$.
		
		By assumption, there exist $s \mathrm e^{\mathrm i \beta} \in 
		W_\varepsilon \cap \{ r \mathrm e^{\mathrm i \alpha}: r \in [0,1],
		\alpha \in [\varphi, \varphi + \vartheta] \}$ and $w \in W$ with
		$|w-s \mathrm e^{\mathrm i \beta}| \leq \varepsilon$. It is easy to
		see that $s \geq 1 - \delta - \varepsilon$. Finally,
		\begin{align*}
			|w - \mathrm e^{\mathrm i \varphi}| &\leq
				|w-s \mathrm e^{\mathrm i \beta}| +
				|s \mathrm e^{\mathrm i \beta} - 
				s \mathrm e^{\mathrm i \varphi}|+
				|s\mathrm e^{\mathrm i \varphi} - 
				\mathrm e ^{\mathrm i \varphi}|\\
			&\leq \varepsilon + \vartheta + (\delta +\varepsilon)
				= 2\varepsilon + \delta + \vartheta. \qedhere
		\end{align*}				
	\end{proof}
\end{lem}
\begin{prop}
	\label{ProductOfFiniteGroups}
	Let $A$ be a compact abelian group, let $(G_l)_{l \in \mathbb N}$ be a
	sequence of finite abelian groups, set $G = \prod_{l = 1}^\infty G_l
	\oplus A$, and let $\Lambda$ be an infinite subset of $\widehat G =
	\bigoplus_{l = 1}^\infty \widehat{G_l} \oplus \widehat A$. If 
	$p_{\widehat A}[\Lambda]$ is a finite set, then $T(C_\Lambda(G)) = 2$.
	\begin{proof}
		The beginning is almost like in the proof of Proposition
		\ref{infiniteTorus}.
		
		Fix $f_1, \dotsc, f_n \in S_{C_\Lambda(G)}$ and $\varepsilon > 0$. 
		Since $T_\Lambda(G)$ is dense in $C_\Lambda(G)$, we may assume
		without loss of generality that $f_1, \dotsc, f_n$ are trigonometric
		polynomials. We have to find $g \in S_{C_\Lambda(G)}$ with 
		$\norm{f_k + g}_\infty \geq 2 - \varepsilon$ for $k = 1, \dotsc, n$.
		
		Setting $\Delta = \bigcup_{k = 1}^n \text{spec}f_k$, we get a finite
		subset of $\Lambda$ because every $f_k$ is a trigonometric
		polynomial and therefore has a finite spectrum. Consequently, there
		exists $l_0 \in \mathbb N$ with $p_{\widehat{G_l}}[\Delta] =
		\{ \mathbf{1}_{G_l} \}$	for all	$l > l_0$ and the evaluation of $f_1,
		\dotsc, f_n$ at a point $(x_1, x_2, \dotsc, a) \in G$ just depends on
		the coordinates	$x_1, \dotsc, x_{l_0}$ and $a$.
		
		Since $\widehat{G_1}, \dotsc, \widehat{G_{l_0}}$ and $p_{\widehat A}
		[\Lambda]$ are finite sets and $\Lambda$ is an infinite set, there
		exist an infinite subset $\Lambda_0$ of	$\Lambda$ and elements
		$\gamma_1 \in \widehat{G_1}, \dotsc, \gamma_{l_0} \in
		\widehat{G_{l_0}}, \gamma_A \in \widehat A$ with 
		$p_{\widehat{G_l}}[\Lambda_0] = \{\gamma_l\}$ for $l = 1, \dotsc, l_0$
		and $p_{\widehat A}[\Lambda_0] = \{\gamma_A\}$.	In other words, all
		elements of $\Lambda_0$ coincide in the first $l_0$ coordinates of
		$\bigoplus_{l = 1}^\infty \widehat{G_l}$ and in	the coordinate that
		corresponds to $\widehat A$. We can also assume	that $\Lambda_0$ is a
		Sidon set because every infinite subset of $\widehat G$ contains an
		infinite Sidon set \cite{RudinFourierAnalysis}*{Example 5.7.6.(a)}.
		(Recall that $\Lambda_0$ is said to be a \emph{Sidon set}, if there
		exists a constant $C > 0$ such that $\sum_{\gamma \in \Lambda_0} 
		|\hat f (\gamma)| \leq C \norm{f}_\infty$ for all $f \in
		T_{\Lambda_0}(G)$.) So if $\{ \lambda_1, \lambda_2, \dotsc \}$ is an
		enumeration of $\Lambda_0$, then $(\lambda_s)_{s \in \mathbb N}$ is 
		equivalent to the canonical basis of $\ell^1$.
		
		Set $\gamma = (\overline{\gamma_1}, \dotsc, 
		\overline{\gamma_{l_0}}, \mathbf{1}_{G_{l_0 + 1}}, 
		\mathbf{1}_{G_{l_0 + 2}}, \dotsc, \overline{\gamma_A}) \in 
		\widehat G$. The sequence $(\gamma \lambda_s)_{s \in \mathbb N}$ is
		still equivalent to the	canonical basis of $\ell^1$ and we have for
		every character	$\gamma \lambda_s$ that $p_{\widehat A}
		(\gamma\lambda_s) =	\mathbf{1}_A$ and $p_{\widehat G_l}
		(\gamma\lambda_s) = \mathbf{1}_{G_l}$ for $l = 1, \dotsc, l_0$. Thus
		the evaluation of $\gamma \lambda_1, \gamma \lambda_2, \dotsc$ at a
		point $(x_1, x_2, \dotsc, a) \in G$ does not depend on the coordinates
		$x_1, \dotsc, x_{l_0}$ and $a$.
		
		Choose $n_0 \in \mathbb N$ with $\frac{2\pi}{n_0} \leq 
		\frac \varepsilon 3$ and $\delta \in (0,1)$ with $4n_0 \sqrt \delta
		\leq \frac \varepsilon 3$. By James's $\ell^1$ distortion theorem 
		\cite{AlbiacKaltonTopics}*{Theorem 10.3.1}, there is a normalized
		block basis sequence $(g_s)_{s \in \mathbb N}$ of 
		$(\gamma \lambda_s)_{s \in \mathbb N}$ with
		\begin{equation*}
			(1 - \delta) \sum_{s = 1}^\infty |z_s| \leq
			\norm{\sum_{s = 1}^\infty z_s g_s}_\infty \leq
			\sum_{s = 1}^\infty |z_s|
		\end{equation*}
		for any sequence of complex numbers $(z_s)_{s \in \mathbb N}$ with
		finite support. It follows that for every $n_0$-tuple
		$(z_1, \dotsc, z_{n_0}) \in \mathbb T^{n_0}$ there is $x \in G$ with
		\begin{equation*}
			\left| \sum_{s = 1}^{n_0} z_s g_s(x) \right| \geq n_0 (1 - \delta).
		\end{equation*}
		Using Lemma \ref{Lemma1}, we have for $s,t = 1, \dotsc, n_0$
		\begin{equation*}
			|g_s(x)| \geq 1 - n_0 \delta \quad \text{and} \quad
			|z_s g_s(x) - z_t g_t(x)| \leq 2n_0 \sqrt \delta.
		\end{equation*}
		Setting for $s = 1, \dotsc, n_0$
		\begin{align*}
			W_s &= g_s[G] \cap \{ z \in \mathbb C : |z| \geq 1 - n_0 \delta \}
			\shortintertext{and}
			\widetilde{W_s} &=\{ z \in \mathbb C : \text{there exists 
				$w \in W_s $ with $|w-z| \leq 2n_0 \sqrt\delta$} \},
		\end{align*}
		we conclude that for every tuple $(z_1, \dotsc, z_{n_0}) \in 
		\mathbb T^{n_0}$
		\begin{equation*}
			\bigcap_{s = 1}^{n_0} z_s \widetilde{W_s} \neq \emptyset.
		\end{equation*}
		By Lemma \ref{Lemma2}, there is $s_0 \in \{1, \dotsc, n_0\}$ such that
		for any $\varphi \in [0,2\pi]$
		\begin{equation*}
			\widetilde{W_{s_0}} \cap
			\left\{ r \mathrm e^{\mathrm i \alpha}: r \in [0,1], \alpha \in
			\left[ \varphi,\varphi + \frac{2\pi}{n_0} \right] \right\}
			\neq \emptyset.
		\end{equation*}
		It follows from Lemma \ref{Lemma3} and our choice of $n_0$ and
		$\delta$ that $W_{s_0}$ is an $\varepsilon$-net for $\mathbb T$.
		
		The function $g = \overline{\gamma}g_{s_0}$ is by construction a 
		normalized trigonometric polynomial with spectrum contained in 
		$\Lambda$. Fix $k \in \{1, \dotsc, n\}$. There exists 
		$x^{(k)} = (x^{(k)}_1, x^{(k)}_2, \dotsc, a^{(k)}) \in G$ with
		$|f_k(x^{(k)})| = 1$. By our choice of $g_{s_0}$ we can find 
		$y^{(k)} = (y^{(k)}_1, y^{(k)}_2, \dotsc, b^{(k)}) \in G$ with
		\begin{equation*}
			\left| \gamma(x^{(k)})f_k(x^{(k)})
			- g_{s_0}(y^{(k)}) \right| \leq \varepsilon.
		\end{equation*} 
		Note that $\gamma(x^{(k)})f_k(x^{(k)}) \in \mathbb T$ since $\gamma$ is
		a character. We therefore get
		\begin{align*}
			\norm{f_k + g}_\infty &= \norm{\gamma f_k + g_{s_0}}_\infty\\
			&\geq \left| (\gamma f_k + g_{s_0})
				(x^{(k)}_1, \dotsc, x^{(k)}_{l_0}, y^{(k)}_{l_0+1},
				y^{(k)}_{l_0+2}, \dotsc, a^{(k)}) \right|\\
			&= \left| \gamma(x^{(k)})f_k(x^{(k)})
				+ g_{s_0}(y^{(k)}) \right|\\
			&\geq 2\left| \gamma(x^{(k)})f_k(x^{(k)}) \right| -
				\left| \gamma(x^{(k)})f_k(x^{(k)}) - g_{s_0}(y^{(k)}) \right|\\
			&\geq 2 - \varepsilon. \qedhere			
		\end{align*}
	\end{proof}
\end{prop}
\begin{lem}
	\label{Lemma4}
	Let $G$ be a compact abelian group and let $\gamma \in \widehat G$.
	\begin{enumerate}[\upshape(a)]
		\item If $o(\gamma) = m$, then $\gamma[G] =
			\{ \mathrm e^{2\pi \mathrm i \frac k m}: k = 0, \dotsc, m-1 \}$,
			i.e., the image of $G$ under $\gamma$ is the set of the $m$th 
			roots of unity.
		\item If $o(\gamma) = \infty$, then $\gamma[G] = \mathbb T$.
	\end{enumerate}
	\begin{proof}
		If $o(\gamma) = m$, we have $\gamma(x)^m = 1$ for every $x \in G$. Thus
		every element of $\gamma[G]$ is an $m$th root of unity. Setting 
		$n = |\gamma[G]|$, it follows from Lagrange's theorem that 
		$\gamma(x)^n = 1$ for every $x \in G$. Therefore $n \geq m$ and 
		$\gamma[G]$ has to coincide with $\{ \mathrm e^{2\pi \mathrm i \frac k
		m}: k = 0, \dotsc, m - 1 \}$.
		
		The set $\gamma[G]$ is a compact and therefore closed subgroup of
		$\mathbb T$. Since all proper closed subgroups of $\mathbb T$ are 
		finite \cite{MorrisPontryaginDualityStructureLCAGroups}*{Corollary 3,
		p. 28},	we have $\gamma[G] = \mathbb T$, if $o(\gamma) = \infty$.
	\end{proof}
\end{lem}
\begin{theo}
	\label{GeneralTheorem}
	Let $G$ be a compact abelian group and let $\Lambda$ be an infinite subset
	of $\widehat G$. Then $T(C_\Lambda(G)) = 2$.
	\begin{proof}
		We start like in the proofs of Proposition \ref{infiniteTorus} and
		\ref{ProductOfFiniteGroups}.
		
		Fix $f_1, \dotsc, f_n \in S_{C_\Lambda(G)}$ and $\varepsilon > 0$. 
		Since $T_\Lambda(G)$ is dense in $C_\Lambda(G)$, we may assume
		without loss of	generality that $f_1, \dotsc, f_n$ are trigonometric
		polynomials. We have to find $g \in S_{C_\Lambda(G)}$ with 
		$\norm{f_k + g}_\infty \geq 2 - \varepsilon$ for $k = 1, \dotsc, n$.
		
		Setting $\Delta = \bigcup_{k = 1}^n \text{spec}f_k$, we get a finite
		subset of $\Lambda$ because every $f_k$ is a trigonometric
		polynomial and therefore has a finite spectrum.
		
		We can assume, by passing to a countably infinite subset if
		necessary, that $\Lambda$ is countable. Hence $\Gamma =	\langle
		\Lambda \rangle$, the group generated by $\Lambda$, is a countable
		subgroup of	$\widehat G$.
		
		Let $M$ be a maximal independent subset of $\Gamma$ and let $
		\Gamma_1 = \langle M \rangle$ be the subgroup of $\Gamma$ that is 
		generated by $M$. Defining inductively $\Gamma_l =
		\{ \gamma \in \Gamma: \gamma^l \in \Gamma_{l - 1}\}$ for $l = 2, 3,
		\dotsc$, we get an increasing sequence $(\Gamma_l)_{l \in \mathbb N}$
		of subgroups of $\Gamma$. Since $M$ is a maximal independent subset of
		$\Gamma$, we have that $\bigcup_{l = 1}^\infty \Gamma_l = \Gamma$. 
		Furthermore, every $\Gamma_l$ is a direct sum of cyclic groups
		\cite{FuchsInfiniteAbelianGroupsI}*{Corollary 18.4}. We distinguish
		two cases depending on whether or not there exists $\Gamma_l$ that
		contains $\Delta$ and infinitely many elements of $\Lambda$.
		
		\emph{First case:} Suppose that there exists $l_0 \in \mathbb N$ such
		that $\Delta \subset \Gamma_{l_0}$ and $\Lambda_0 = \Gamma_{l_0}\cap
		\Lambda$ is an infinite set.
		
		By our choice of $\Gamma_{l_0}$, the functions $f_1, \dotsc, f_n$ and
		all characters $\gamma \in \Lambda_0$ are constant on the cosets
		of $G/(\Gamma_{l_0})^\perp$ and can therefore be considered as
		functions and characters on $G_0 = G/(\Gamma_{l_0})^\perp$. (To
		simplify notation, we continue to write $f_1, \dotsc, f_n$.) Note that
		$\Gamma_{l_0}$ is the dual group of $G_0$. Since $\Gamma_{l_0}$ is
		a direct sum of cyclic groups, there exists a sequence
		$(\widehat{G_s})_{s \in \mathbb N}$ of finite abelian groups such that
		$\Gamma_{l_0} =	\mathbb Z^{(\mathbb N)} \oplus 
		\bigoplus_{s = 1}^\infty \widehat{G_s}$	or  $\Gamma_{l_0} = 
		\mathbb Z^{n_0} \oplus \bigoplus_{s = 1}^\infty	\widehat{G_s}$ for 
		adequate $n_0 \in \mathbb N$. Hence	$G_0 = \mathbb T^{\mathbb N} 
		\oplus \prod_{s = 1}^\infty G_s$ or  $G_0 = \mathbb T^{n_0} \oplus
		\prod_{s = 1}^\infty G_s$. Let $p_1, p_2, \dotsc$ be the projections 
		from $\Gamma_{l_0}$ onto $\mathbb Z$.
		
		If there exists $s_0 \in \mathbb N$ such that $p_{s_0}[\Lambda_0]$ 
		contains infinitely many elements or if there exist arbitrarily large
		$s \in \mathbb N$ with $p_s[\Lambda_0] \neq \{0\}$, then
		$T(C_{\Lambda_0}(G_0)) = 2$ by Proposition \ref{finiteTorus} or
		\ref{infiniteTorus}. Otherwise $p_{\mathbb Z ^{(\mathbb N)}}
		[\Lambda_0]$ (or $p_{\mathbb Z ^{n_0}}[\Lambda_0]$) is a finite set and
		$T(C_{\Lambda_0}(G_0)) = 2$ by Proposition \ref{ProductOfFiniteGroups}.
		So we can find $\tilde g \in S_{C_{\Lambda_0}(G_0)}$ with
		$\norm{f_k + \tilde g}_\infty \geq 2 - \varepsilon$ for $k = 1,
		\dotsc, n$. Setting $g=\tilde g \circ \pi$ where
		$\pi$ is the canonical map from $G$ onto $G_0 = G/(\Gamma_0)^\perp$, we
		get $\norm{f_k + g}_\infty \geq 2 - \varepsilon$ for $k = 1,\dotsc,
		n$.
		
		\emph{Second case:} Suppose that there exist arbitrarily large
		$l \in \mathbb N$ with $(\Gamma_l \setminus \Gamma_{l - 1}) \cap
		\Lambda \neq \emptyset$.
		
		Fix $l_0 \in \mathbb N$ with $\Delta \subset \Gamma_{l_0}$ and choose
		$l_1 \in \mathbb N$ with $l_1 > l_0^2$, $\frac{2\pi}{l_1} \leq 
		\varepsilon$ and $(\Gamma_{l_1} \setminus \Gamma_{l_1 - 1}) \cap
		\Lambda \neq \emptyset$. By our choice of $\Gamma_{l_0}$, the
		functions $f_1, \dotsc, f_n$ are constant on the cosets of
		$G/(\Gamma_{l_0})^\perp$ and therefore
		\begin{equation}
			\label{GeneralTheoremEquation1}
			f_k(xy) = f_k(x) \quad (x \in G, y \in (\Gamma_{l_0})^\perp)
		\end{equation}
		for $k = 1, \dotsc, n$. Pick $g \in (\Gamma_l \setminus 
		\Gamma_{l - 1}) \cap \Lambda$ and denote by $\tilde g$ the restriction 
		of $g$ to $(\Gamma_{l_0})^\perp$. What can we say about the order of 
		$\tilde g$? Since $(\Gamma_{l_0})^{\perp \perp} = \Gamma_{l_0}$, we 
		have for every $m \in \mathbb N$ that $\tilde g ^m = 
		\mathbf 1_{(\Gamma_{l_0})^\perp}$ if and only if 
		$g^m \in \Gamma_{l_0}$.
		
		Suppose that $\tilde g ^m = \mathbf 1_{(\Gamma_{l_0})^\perp}$ for
		some $2 \leq m \leq l_0$. Then $\tilde g ^{ml_0} = 
		\mathbf 1_{(\Gamma_{l_0})^\perp}$ as well and 
		$g^{ml_0} \in \Gamma_{l_0}$.
		Consequently, $g \in \Gamma_{ml_0}$ because $g^{ml_0} \in \Gamma_{l_0}
		\subset \Gamma_{ml_0 - 1}$. But this contradicts our choice of $g$ and
		$l_1$ because $l_1 > ml_0$.	Assuming that $\tilde g^m = 
		\mathbf 1_{(\Gamma_{l_0})^\perp}$ for some $l_0 < m < l_1$ leads to the same
		contradiction. The order of $\tilde g$ is therefore at least $l_1$.
		By our choice of $l_1$ and by Lemma \ref{Lemma4}, we get that
		$\tilde g[(\Gamma_{l_0})^\perp]$ is an $\varepsilon$-net for $\mathbb
		T$.
		
		Fix now $k \in \{1, \dotsc, n\}$ and choose $x^{(k)} \in G$ with
		$|f_k(x^{(k)})| = 1$ and $y^{(k)} \in (\Gamma_{l_0})^\perp$ with
		\begin{equation}
			\label{GeneralTheoremEquation2}
			\left| f_k(x^{(k)}) - g(x^{(k)})\tilde g (y^{(k)}) \right|
			\leq \varepsilon.
		\end{equation}
		Note that $g$ is a character and hence $g(x^{(k)}) \in \mathbb T$.
		Using (\ref{GeneralTheoremEquation1}) and
		(\ref{GeneralTheoremEquation2}), we get
		\begin{align*}
			\norm{f_k + g}_\infty &\geq \left| f_k(x^{(k)}y^{(k)}) +
				g(x^{(k)}y^{(k)}) \right|\\
			&= \left| f_k(x^{(k)}) + g(x^{(k)})\tilde g (y^{(k)}) \right|\\
			&\geq 2 \left| f_k(x^{(k)}) \right| - \left| f_k(x^{(k)}) -
				g(x^{(k)})\tilde g (y^{(k)}) \right|\\
			&\geq 2 - \varepsilon. \qedhere		
		\end{align*}		
	\end{proof}
\end{theo}
\begin{cor}
	Let $G$ be a metrizable, compact abelian group and let $\Lambda$ be a
	subset of $\widehat G$. The space $C_\Lambda(G)$ has the almost Daugavet
	property if and only if $\Lambda$ contains infinitely many elements.
	\begin{proof}
		Every almost Daugavet space is infinite-dimensional and so the 
		condition is necessary.
		
		If $G$ is a metrizable, compact abelian group, then $\widehat G$ is
		countable \cite{RudinFourierAnalysis}*{Theorem 2.2.6} and $C(G)$ is 
		separable. Since for separable Banach spaces the almost Daugavet
		property can be characterized via the thickness
		\cite{KadetsShepelskaWernerThicknessAlmostDaugavet}*{Theorem 1.1}, 
		it is sufficient to prove that $T(C_\Lambda(G)) = 2$. But this is
		given by Theorem \ref{GeneralTheorem}.			 
	\end{proof}
\end{cor}

\section{Subspaces of $L$-embedded spaces}

To deal with translation-invariant subspaces of $L^1(G)$ we will consider a
more general class of Banach spaces. A linear projection $P$ on a Banach space
$X$ is called an \emph{$L$-projection}, if
\begin{equation*}
	\norm{x}=\norm{Px} + \norm{x-Px} \quad (x \in X).
\end{equation*} 
A closed subspace of $X$ is called an \emph{$L$-summand}, if it is the range of
an $L$-projection, and $X$ is called \emph{$L$-embedded}, if $X$ is an 
$L$-summand in $X^{**}$. Classical examples of $L$-embedded spaces are
$L^1(\mu)$-spaces, preduals of von Neumann algebras, and the Hardy space $H^1$
\cite{HarmandWernerMideale}*{Example IV.1.1}.

Using the principle of local reflexivity, it is easy to see that a
non-reflexive, $L$-embedded space has thickness 2. We will strengthen this
and will show that every non-reflexive subspace of an $L$-embedded space has 
thickness 2. Let us recall the following result from the theory of
$L$-embedded spaces \cite{HarmandWernerMideale}*{claim in the proof of Theorem
IV.2.7}}
\begin{prop}
	\label{PropositionNearAccumulationPoint}
	Let $X$ be an $L$-embedded space with $X^{**} = X\oplus_1 X^s$, 
	let $\varepsilon$
	be a number with $0 < \varepsilon < \frac 1 4$, and let 
	$(y_l)_{l \in \mathbb N}$ be a sequence in $X$ with
	\begin{equation*}
		(1 - \varepsilon)\sum_{l = 1}^\infty |a_l| \leq 
		\norm{\sum_{l = 1}^\infty a_l y_l} \leq
		\sum_{l = 1}^\infty |a_l|
	\end{equation*}
	for any sequence of scalars $(a_l)_{l \in \mathbb N}$ with
	finite support. Then there exists $x_s \in X^s$ such that
	\begin{equation*}
		1 - 4\sqrt \varepsilon \leq \norm{x_s} \leq 1
	\end{equation*}
	and for all $\delta > 0$, all $x_1^*, \dotsc, x_n^* \in X^{**}$ and all
	$l_0 \in \mathbb N $ there is $l \geq l_0$ with
	\begin{equation*}
	|x_s(x_k^*) - x_k^*(y_l)| \leq 3 \sqrt \varepsilon \norm{x_k^*} + \delta
	\quad (k = 1, \dotsc, n).	
	\end{equation*}
	In other words, there is $x_s \in X^s$ which is ``close'' to a weak*
	accumulation point of $(y_l)_{l \in \mathbb N}$.
\end{prop}
\begin{theo}
	\label{LEmbedded}
	Let $X$ be an $L$-embedded space with $X^{**} = X \oplus_1 X^s$ and let 
	$Y$ be a closed subspace of $X$ which is not reflexive. Then $T(Y) = 2$.
	\begin{proof}
		Fix $x_1, \dotsc, x_n \in S_Y$ and $\varepsilon > 0$. We have to find 
		$y \in S_Y$ with $\norm{x_k + y} \geq 2 - \varepsilon$ for 
		$k = 1, \dotsc, n$.
		
		Choose $\delta > 0$ with $7 \sqrt \delta + 2 \delta \leq \varepsilon$.
		Every non-reflexive subspace of $X$ contains a copy of $\ell^1$
		\cite{HarmandWernerMideale}*{Corollary IV.2.3} and by James's $\ell_1$
		distortion theorem \cite{AlbiacKaltonTopics}*{Theorem 10.3.1} there is
		a sequence $(y_l)_{l \in \mathbb N}$ in $Y$ with
		\begin{equation*}
			(1 - \delta)\sum_{l = 1}^\infty |a_l| \leq 
			\norm{\sum_{l = 1}^\infty a_l y_l} \leq
			\sum_{l = 1}^\infty |a_l|
		\end{equation*}
		for any sequence of scalars $(a_l)_{l \in \mathbb N}$ with
		finite support. Let $x_s \in X^s$ be ``close'' to a weak* accumulation
		point of $(y_l)_{l \in \mathbb N}$ as in Proposition
		\ref{PropositionNearAccumulationPoint}. Since $X^{**} = X \oplus_1
		X^s$, we have for $k = 1, \dotsc, n$
		\begin{equation*}
			\norm{x_k + x_s} = \norm{x_k} + \norm{x_s} \geq 2 - 4\sqrt \delta.
		\end{equation*}
		Thus there exist functionals $x_1^*, \dotsc, x_n^* \in S_{X^*}$ with
		\begin{align*}
			|x_k^*(x_k) + x_s(x_k^*)| &\geq 2 - 4 \sqrt \delta - \delta
			\shortintertext{and $l \in \mathbb N$ with}
			|x_s(x_k^*) - x_k^*(y_l)| &\leq 3 \sqrt \delta + \delta
		\end{align*}
		for $k = 1,\dotsc, n$.
		
		Fix $k \in \{1, \dotsc, n\}$. Using the last two inequalities leads to
		\begin{align*}
			\norm{x_k + y_l} &\geq |x_k^*(x_k) + x_k^*(y_l)|\\
			&\geq |x_k^*(x_k) + x_s(x_k^*)| - |x_s(x_k^*) - x_k^*(y_l)|\\
			&\geq (2 - 4\sqrt \delta - \delta) - (3 \sqrt \delta + \delta)\\
			&\geq 2 - \varepsilon. \qedhere  
		\end{align*}
	\end{proof}
\end{theo}
\begin{cor}
	Let $X$ be an $L$-embedded space and let $Y$ be a separable, closed
	subspace of $X$. If $Y$ is not reflexive, then $Y$ has the almost
	Daugavet property.
	\begin{proof}
		The space $Y$ has thickness 2 by Theorem \ref{LEmbedded} and this is
		for separable spaces equivalent	to the almost Daugavet property
		\cite{KadetsShepelskaWernerThicknessAlmostDaugavet}*{Theorem 1.1}.
	\end{proof}
\end{cor}
Let us use this result in the setting of translation-invariant subspaces of
$L^1(G)$. Suppose that $G$ is a compact abelian group, $\Lambda$ a subset of
its dual group $\widehat G$ and $0 < r < p < \infty$. The set $\Lambda$ is
said to \emph{be of type} $(r,p)$, if there is a constant $C > 0$ such that
\begin{equation*}
	\norm{f}_p \leq C \norm{f}_r
\end{equation*}
for every $f \in T_\Lambda(G)$. In other words, if $\norm{\,\cdot\,}_r$ and
$\norm{\,\cdot\,}_p$ are equivalent on $T_\Lambda(G)$. We say furthermore that
$\Lambda$ is a \emph{$\Lambda(p)$ set}, if $\Lambda$ is of type $(r,p)$ for
some $r < p$.
\begin{cor}
	Let $G$ be a metrizable, compact abelian group and let $\Lambda$ be a
	subset of $\widehat G$. The space $L^1_\Lambda(G)$ has the almost Daugavet
	property if and only if $\Lambda$ is not a $\Lambda(1)$ set.
	\begin{proof}
		Every almost Daugavet space contains a copy of $\ell^1$
		\cite{KadetsShepelskaWernerThicknessAlmostDaugavet}*{Corollary 3.3} and
		is therefore not reflexive. So the condition is necessary because
		$L^1_\Lambda(G)$ is reflexive if and only if $\Lambda$ is a
		$\Lambda(1)$ set \cite{HareElementaryProofLambdap}*{Corollary}.
		
		If $G$ is a metrizable, compact abelian group, then $\widehat G$ is
		countable \cite{RudinFourierAnalysis}*{Theorem 2.2.6} and 
		$L^1(G)$ is separable. If $\Lambda$ is not a $\Lambda(1)$ set, then
		$L^1_\Lambda(G)$ is not reflexive and $T(L^1_\Lambda(G)) = 2$
		by Theorem \ref{LEmbedded}. But this is for separable spaces equivalent
		to the almost Daugavet property
		\cite{KadetsShepelskaWernerThicknessAlmostDaugavet}*{Theorem 1.1}.
	\end{proof}	
\end{cor}
\section{Remarks}
The almost Daugavet property is strictly weaker than the Daugavet
property for translation-invariant subspaces of $C(G)$ or $L^1(G)$. If we
set $\Lambda = \{ 3^n : n \in \mathbb N \}$, then $\Lambda$ is a Sidon set. So
$C_\Lambda(\mathbb T)$ is isomorphic to $\ell^1$, has the 
Radon-Nikod\'ym property and therefore not the Daugavet property. 
But $\Lambda$ is an
infinite set and $C_\Lambda(\mathbb T)$ has the almost Daugavet property.
Analogously, $L^1_\mathbb N(\mathbb T)$ is isomorphic to the Hardy space
$H^1_0$, has therefore the Radon-Nikod\'ym property and fails the Daugavet
property. But $\mathbb N$ is not a $\Lambda(1)$ set and $L^1_\mathbb N(G)$ has
the almost Daugavet property.

We say that a Banach space $X$ has the fixed point property, if given any 
non-empty, closed, bounded and convex subset $C$ of $X$, every non-expansive
mapping $T: C \rightarrow C$ has a fixed point. Here $T$ is non-expansive,
if $\norm{Tx - Ty}\leq \norm{x - y}$ for all $x,y \in C$. By considering
$C = \{ (x_n)_{n \in \mathbb N} \in S_{\ell^1}: x_n \geq 0 \}$ and the right
shift operator, it can be shown that $\ell^1$ does not have the fixed point
property \cite{DowlingLennardNonreflexiveSubspaceL1FailsFixedPointProperty}*
{Theorem 1.2}. This counterexample can be transferred to all Banach spaces
that contain an \emph{asymptotically isometric copy of $\ell^1$}. A Banach
space $X$ is said to contain an asymptotically isometric copy of $\ell^1$,
if there is a null sequence $(\varepsilon_n)_{n \in \mathbb N}$ in $(0,1)$
and a sequence $(x_n)_{n \in \mathbb N}$ in $X$ such that
\begin{equation*}
	\sum_{n = 1}^\infty (1 - \varepsilon_n)|a_n| \leq
	\norm{\sum_{n = 1}^\infty a_nx_n} \leq \sum_{n = 1}^\infty |a_n|
\end{equation*}
for any sequence of scalars $(a_n)_{n \in \mathbb N}$ with finite support. 
Every Banach space $X$ with $T(X) = 2$ contains an asymptotically isometric
copy of $\ell^1$ \cite{KadetsShepelskaWernerThicknessAlmostDaugavet}*
{implicitly in the proof of Propositions 3.2 and 3.4}. So Theorem
\ref{LEmbedded} gives another proof of the fact that every non-reflexive 
subspace of $L^1[0,1]$ or more generally every non-reflexive subspace of an
$L$-embedded space fails the fixed point property (cf. \citelist{
\cite{DowlingLennardNonreflexiveSubspaceL1FailsFixedPointProperty}*
{Theorem 1.4}\cite{PfitznerNoteAsymptoticallyIsometricCopiesl1c0}*
{Corollary 4}}).

\section*{Acknowledgements}

This is part of the author's Ph.D. thesis, written under the supervision of
D. Werner at the Freie Universit\"at Berlin.

\end{document}